\documentclass[twoside,leqno]{amsart}

\usepackage{amssymb}
\usepackage{amsmath}
\usepackage{graphicx}
\usepackage{textcomp}
\usepackage{xcolor}
\usepackage{hyperref}
\usepackage{stmaryrd}

\newcounter{item}[section]

\renewcommand{\theitem}{\thesection.\arabic{item}}

\newcommand{\qitem}[2]{\refstepcounter{item}
  {\bf \theitem\  #1.}\  {#2 }}

\newcommand{\eitem}{\medskip}
    \newcommand{\bef}[4]{\begin{figure}[#4]      % 1== caption  2== vspace
\refstepcounter{figure}
\label{#3}                                   % 3== name of  file
\begin{center}

\includegraphics[height=#2]{#3}
\medskip

 {\normalsize\sl Figure \thefigure. #1   }

\end{center}

  }  %the "s" in [s] can be replaced

    \newcommand{\ef}{\end{figure}}

\newcommand{\cI}{\mathcal{I}}
\newcommand{\eps}{\varepsilon}

\newcommand{\R}{\mathbb{R}}

\newcommand{\E}{{\rm E\, }} %expectation

\newcommand{\dsp}{\displaystyle}
\newcommand{\be}{\begin{equation} }

\newcommand{\ee}{\end{equation}}
\newcommand{\beq}{\begin{eqnarray} }
\newcommand{\eeq}{\end{eqnarray}}
\newcommand{\beqn}{\begin{eqnarray*} }
\newcommand{\eeqn}{\end{eqnarray*}}

\newcommand{\di}{{\: \rm d}}   % This is "d" in math formulas with extra sp.
\newcommand{\fai}{\raisebox{.1em}{$\varphi$}}

\newcommand{\limitsi}[2]{
\left |_{\raisebox{-.5em}{\scriptsize $#1$}}
^{\raisebox{.5em}{\scriptsize $#2$}} \rule{0em}{1em}
\right .
}

\newcommand{\expd}{{\rm Exp\, }} %exponential distribution

\begin{document}

\begin{center}

{\LARGE\bf A tight   Hermite-Hadamard's inequality  and a generic method for comparison between  residuals of inequalities with  convex functions}

\vspace{3em}

{\sc Milan Merkle and  Zoran D. Mitrovi\'c}

\bigskip

\parbox{25cc}{
{\bf Abstract. }{\small  We present a  tight parametrical  Hermite-Hadamard type  inequality with probability measure, which yields a
considerably closer
upper bound for the mean value of  convex function than the classical one.
Our inequality becomes equality not only with affine functions, but also with a family of V-shaped curves determined
by the parameter. The residual (error) of this inequality is strictly smaller than in the classical Hermite-Hadamard inequality
 under any probability measure and with all non-affine convex functions. In the framework of  Karamata's theorem on the inequalities with
convex functions, we propose a method of measuring a global performance of inequalities in terms of average residuals over functions of the type
 $x\mapsto |x-u|$.  Using average residuals  enables comparing  two or more inequalities as themselves,
 with same or different measures and  without referring to a particular function.
Our method is applicable to all Karamata's type  inequalities, with integrals or sums.
A numerical experiment  with three different measures indicates that the average residual in our inequality is about 4 times smaller than
in classical right Hermite-Hadamard, and also is smaller than in Jensen's inequality, with all three measures.

\medskip


\medskip


} }

\end{center}

\section{Introduction}\label{intro}

\setcounter{section}{1}

For a  non-negative measure  $\mu$ on $[a,b]$,  such that $\mu [a,b]=1$ (probability measure), let $c=\int x\di \mu(x)$.
From long ago \cite{luppef76, McShane37}, it is known that for a convex function $f$ it holds

\be
\label{genhh}
 f(c)\leq \int f(x)\di \mu (x) \leq \frac{b-c}{b-a}f(a)+ \frac{c-a}{b-a}f(b).
 \ee

A  case with $\mu$ being the Lebesgue probability measure on $[a,b]$, with $\di \mu(x) = \frac{1}{b-a}\di x$,
is  originally stated by C. Hermite  and J. S. Hadamard independently in late 19th century (see \cite{mitlac} for more history):

\be
\label{hed}
f\left(\frac{a+b}{2}\right) \leq\frac{1}{b-a} \int_a^b f(x)\di x \leq \frac{f(a)+f(b)}{2}.
\ee

Hermite-Hadamard ($HH$) inequality from the beginning has been used in problems of approximations the integral in the middle, using
left inequality (midpoint rule) or  the right one (trapezoid rule). It is well known that the residual (error) in the right inequality
is larger than in the left one, see \cite{bull78,hamm58}, and  there is a voluminous literature on  refinement of the right side of
(\ref{hed}) like in  \cite{allas2015,guess2002,guess2018,mm99,merc2009,Olbs2015} and many more. Regardless of applications, it is always
desirable to have an inequality with smaller residual. This paper offers two contributions to this topic. In  Section 2 we present a new parametrical right
 bound in (\ref{genhh}), which gives much smaller residual
for all measures and all non-affine convex functions, with all values of the parameter and without any additional assumptions. In Section 3, we develop a
method via Karamata's theorem (know also as Levin-Ste\v ckin) for estimating the
residuals of inequalities for convex functions, and comparing residuals of different inequalities. To the best of our  knowledge, this is the unique
method capable to compare any two or more inequalities globally, without referring to a certain function.  Numerical experiments in Section 4 confirm
the theoretical results and also indicate that the residual in our  tight  inequality is smaller not only in comparison with the right bound in \ref{genhh}
but also with respect to the left bound.

 In this paper we adopt the setup  with countably additive probability measures on $\R$ and
 Lebesgue-Stieltjes integrals (as the most general integral that incorporates Riemann and Riemann-Stieltjes one) over a compact interval $[a,b]$.
 To avoid repetitions, let us state several notions and conditions.

\qitem{Notions and conditions}{\label{notco} In the rest of this paper, by  measure will be understood a countably additive probability measure
$\mu$ on Borel sigma algebra on $\R$ such that $\mu [a,b]=1$, $a<b$. Let $X$ be the random variable  associated to $\mu$, with distribution
function $G$ defined as $G(x)= \mu (-\infty,x]$. The integral of an integrable function $f$
with respect to measure $\mu$ is expressed as the Lebesgue-Stieltjes integral  $\int f(x)\di G(x)$, or in  more compact terms
of expectation operator and random variables, as $\E (f(X))$.  Under "convex function on $[a,b]$" we understand a function which is convex
on some open interval $I$ that contains $[a,b]$. \hfill$\square$
}

The left inequality  in (\ref{genhh}) is Jensen's inequality originally proved by Jensen \cite{jen06} and generalized by McShane \cite{McShane37}.
A very simple proof can be find in \cite{durp}. The following theorem presents Jensen's inequality on the compact interval, in our setup in \ref{notco}.

\eitem

\qitem{Theorem}{\label{tjens} Let $\mu$ and $G$ be as in \ref{notco}, and let $f$ be a convex function on  $[a,b]$. Then

\be
\label{jens}% J
f\left( \int_{[a,b]} x\di G(x)\right) \leq \int_{[a,b]} f(x)\di G(x), \quad \mbox{or equivalently,}\quad   f(\E X) \leq \E (f(X))
\ee

}
\eitem

The right inequality in (\ref{genhh}) is generalized by A. Lupa\c{s}  \cite{luppef76} with an abstract linear functional.
Here we give a  formulation under our setup, and a short direct proof.

\qitem{Theorem}{\label{thad} Under the same conditions as in Theorem \ref{tjens}, for every  convex function $f$ it holds:
\be
\label{had} %H
\int_{[a,b]} f(x)\di G(x)  \leq \frac{b-\int_{[a,b]}x \di G(x)}{b-a}f(a) + \frac{\int_{[a,b]} x \di G(x)-a}{b-a}f(b),
\ee
or equivalently,

\[ \E f(X) \leq \frac{b-\E X}{b-a}f(a) + \frac{\E X-a}{b-a}f(b)\]

 Moreover, if
\be
\label{hadc}
 \E f(X)\leq \alpha f(a) + \beta f(b)
 \ee
 holds for every  convex function $f$, then  $\beta=1-\alpha$ and $\alpha a +\beta b=\E (X)$.

 \begin{proof} For a convex function $f$ on $[a,b]$, it holds
\[ f(x) \leq \frac{b-x}{b-a}f(a) + \frac{x-a}{b-a} f(b), \]
and after the integration on both sides, we get (\ref{had}).  The second part follows by plugging  $f=x$, $f=-x$ and $f=1$.
\end{proof}

 }
\eitem

In what follows we  use abbreviation $J$ for the Jensen's inequality (\ref{jens}) and $H$ for the one in (\ref{had}). A complete double-side
inequality will be denoted as $HH$. In the next section  we  present a
$H$-type inequality with the right term being  closer to $\E f(X)$ than in (\ref{had}).
In Section 3   we  use the Karamata's theorem  \cite{karam32}, see also \cite[page 645]{marolk} to  define a method  that can be used for
comparison among inequalities with convex functions. In Section 4 we show the results of some numerical experiments and comparisons.

\section{A new tight H-type inequality\label{tesno}}

From Theorem \ref{thad} it follows that the  $H$-inequality can not be improved by changing the weights associated to $f(a)$ and $f(b)$. However,
if we add an arbitrary point $t\in (a,b)$ and re-calculate the weights, the  sum on the right hand side becomes considerably closer to $\E f(X)$ than in
$H$-inequality, for all underlying measures and all  non-affine convex functions.

\smallskip

\qitem{Assumption\label{ass}}{ In order to avoid separations of cases,
in this section we exclude measures concentrated on less than 3 points in $[a,b]$, i.e., we assume that there are no
$x_1,x_2\in [a,b]$ such that $\mu \{x_1\}+\mu \{x_2\}=1$ and $\mu \{x_1\}\geq 0$, $\mu \{x_2\}\geq 0$.
}

\eitem

\qitem{Theorem}{\label{thh}  Let $f$ be a convex function on a compact interval $[a,b]$, $a<b$, and  let
 $\mu$ be a probability measure with distribution function $G$, under notations and conditions as in \ref{notco}. Then the following holds:

(i) For any fixed $t\in (a,b)$,

\beq
\label{thh1}
\int_{[a,b]} f(x)\di G(x) & \leq &  \frac{f(a)}{t-a}\int_{[a,t]} (t-x)\di G(x) + \frac{f(b)}{b-t}\int_{(t,b]}(x-t)\di G(x)\nonumber\\
                         & & + f(t)\left( \frac{1}{t-a}\int_{[a,t]} (x-a)\di G(x) + \frac{1}{b-t}\int_{(t,b]}(b-x)\di G(x).   \right)
\eeq

(ii) Equivalently, for fixed $t$ let  $\lambda =\lambda(t)= \frac{b -t}{b-a}$, so that  $t= \lambda a +(1-\lambda) b$. Then

\beq
\label{thh2}
\int_{[a,b]} f(x)\di G(x) & \leq &\frac{b-\int x\di G(x)}{b-a} f(a) + \frac{\int x\di G(x)-a}{b-a} f(b) \nonumber\\
                          &      & +\left(f(t)-\lambda f(a)-(1-\lambda)f(b)\right)\\
                          &      &\times \left(\frac{1}{t-a}\int_{[a,t]}(x-a)\di G(x) +\frac{1}{b-t} \int_{(t,b]} (b-x)\di G(x)\right) \nonumber
\eeq

\begin{proof}
The convexity of  $f$
implies that, for arbitrary $x\in [a,b]$,

\beq
\label{eqcon1}
 f(x) &\leq & \frac{t -x}{t-a}f(a) + \frac{x-a}{t-a}f(t),\quad a\leq x\leq t,\\
 \label{eqcon2}
 f(x) & \leq & \frac{b -x}{b-t}f(t) + \frac{x-t}{b-t}f(b),\quad t\leq x\leq b .
 \eeq
 Integrating (\ref{eqcon1}) with respect to measure $\mu$ on $[a,t]$ and (\ref{eqcon2}) on
 $(t,b]$, and adding, we get

\beq
\label{hadii}
\int_{[a,b]} f(x)\di \mu (x) &\leq & \frac{\E (t -X)\cdot I_{[a,t]}(X)}{t  -a} f(a) +\frac{\E (X - a) \cdot I_{[a,t]}(X)}{t -a} f(t)\nonumber \\
                              & & + \frac{\E (b -X)\cdot  I_{(t,b]}(X)}{b-t} f(t) +\frac{\E( X -t) \cdot I_{(t,b]}(X)}{b-t} f(b),
\eeq

which is the inequality (\ref{thh1}) in terms of random variables. To show  the equivalence between (\ref{thh1}) and (\ref{thh2}), it suffices to
verify that coefficients with $f(a), f(b)$ and $f(t)$ are equal in both formulae.
\end{proof}
}
\eitem

 Unlike Jensen's inequality, the right-hand side of $HH$
inequality has not been much used so far in probability and statistics. Nevertheless, it might be of interest to formulate the new complete
$HH$-type inequality in terms of random variables as a corollary of   theorem \ref{thh}.

\qitem{Corollary}{\label{Colt} Let $X$ be a random variable supported on $[a,b]$, $a<b$ and with distribution function $G$. For any $t\in (a,b)$ and
a convex function $f$, it holds

\[ f(\E X) \leq \E f(X) \leq p_a f(a) +p_t f(t) + p_b f(t), \]
where

\[p_a=\frac{\int_{[a,t]} (t -x)\di G(x)}{t -a}, \quad p_b = \frac{\int_{(t,b]} ( x -t) \di G(x)}{b-t}\]

\[ p_t= \frac{\int_{[a,t]}(x - a) \di G(x)}{t -a} + \frac{\int_{(t,b]}(b -x)\di G(x)}{b-t}, \]
and  $p_a+p_b+p_c=\int \di G(x) =1$ \hfill$\square$
}
\eitem

The inequality proved in Theorem \ref{thh} will be referred to as  Tight Hermite-Hadamard (abbreviated $TH$) inequality. Let $R_{J}$, $R_{H}$ and
$R_{TH}$ be the corresponding  residuals in $J$, $H$ and $TH$ inequalities. Given the interval $[a,b]$,
the size of residuals depends on the underlying measure $\mu$, and on the function $f$.

\qitem{Lemma}{\label{tlema} For any convex function on $[a,b]$ and  any measure that satisfies assumption \ref{ass}, it holds

{\bf a)} $R_{H} (\mu, f)=0$ if and only if $f$ is affine function. The same holds for $R_J(\mu, f)$.

{\bf b)}  $ R_{TH} (\mu,f,t)< R_{H} (\mu, f)$  for all convex non-affine functions, for all $t\in (a,b)$.

{\bf c)}  $R_{TH} (\mu,f,t)=0$  if and only if
\be
\label{tiff}
f(x)=\left( \frac{t-x}{t-a}\alpha + \frac{x-a}{t-a}\tau\right) I_{[a,t]}(x)
    + \left( \frac{b-x}{b-t}\tau + \frac{x-t}{b-t}\beta\right) I_{(t,b]}(x),
\ee
for some real numbers $\alpha, \beta, \tau$.

\smallskip

Note that the functions defined by (\ref{tiff}) are either affine or their graphs are V-shaped,
with two lines that meet at the point $(t,\tau)$ and with endpoints  $(a,\alpha)$ and $(b,\beta)$.
Such functions are convex if and only if $\tau \leq \min\{\alpha,\beta\}$.

\begin{proof}  We will prove in lemma \ref{compp} that the residuals $R_J$, $R_H$ and $R_{TH}$ are zero for all measures  with an affine $f$,
so we need to prove "only if" part where applies.

{\bf a)} Equality  $R_{H} (\mu, f)=0$  is equivalent to $\E f(X) =p f(a) + (1- p) f(b)$, $p=\frac{b-\E X}{b-a}$. Suppose that $f$ is not
affine. This implies that the graph of $f$ for $x\in (a,b)$ lies  under the chord that connects points $(a,f(a))$ and $(b,f(b))$,  although
the point  $(\E X, \E f(X))$  belongs to  the chord.
This is possible only if $\mu$ is concentrated on the set $\{a,b\}$,
which is excluded by assumption \ref{ass}.  For the proof of necessity for $R_J=0$, see \cite[page 654]{marolk}.

{\bf b)} From  the representation (\ref{thh2}) it follows that

\beq
\label{difere}
\lefteqn{ R_{TH} (\mu,f,t)- R_{H} (\mu, f)  = \left(f(t)-\lambda f(a)-(1-\lambda)f(b)\right)\rule{8cm}{0cm}}\nonumber\\
& &   \times \left(\frac{1}{t-a}\int_{[a,t]}(x-a)\di G(x) +\frac{1}{b-t} \int_{(t,b]} (b-x)\di G(x)\right)
\eeq
According to a), the first term is zero if and only if $f$ is affine; otherwise it is negative. The second
term is  positive under the assumption \ref{ass}, and the claim is proved.

{\bf c)} The function $f$ defined by (\ref{tiff}) satisfies (\ref{eqcon1}) and (\ref{eqcon2}) with equalities. Tracing the proof of Theorem \ref{thh},
the integration with respect to the given measure yields null residual. Moreover, the residual can be zero only if both (\ref{eqcon1}) and (\ref{eqcon2})
are equalities, and this is the case only if $f$ is either affine or in the form (\ref{tiff}).
\end{proof}
}

From Lemma \ref{tlema} it follows that with any measure, inequality $TH$ yields the better approximation to
$\int f(x)\di G(x)$ than the inequality $H$, for every  non-affine convex function.

 \eitem

\qitem{Optimal choice of parameter $t$ for given $f$}{Since inequality $TH$  is valid for any $t\in (a,b)$, it is natural to ask
 which $t$ yields  the smallest residual,  or equivalently, the smallest (negative) difference $ R_{TH} (\mu,f,t)- R_{H} (\mu, f) $
for given $f$ and $\mu$. For given measure $\mu$  with distribution function $G$ and a convex function $f$, this difference can be written
as the function of $\lambda$ using the relation $t=\lambda a + (1-\lambda)b$

\be
\label{ded}
D(\lambda)= R_{TH} (\mu,f,t)- R_{H} (\mu, f)=(f(t)-\lambda f(a)-(1-\lambda)f(b)) \E (g(X)),
\ee
where

\be
\label{deg}
g(x)= \frac{\lambda (x - a)I_{[a,t]}(x)+ (1-\lambda)(b-x)I_{(t,b]}(x)}{\lambda(1-\lambda)(b-a)}, \quad x\in \R,
\ee
and $\lambda = \frac{b-t}{b-a}$. The graph of this function is the continuous triangular curve  which connect points $(a,0)$, $(t,1)$ and $(b,0)$.

A value of  $t$ that minimizes $D(\lambda(t)$ depends on the underlying measure.
In the case of uniform distribution on $[a,b]$, we have $G(x)=\frac{1}{b-a}$ and   $\E (g(X))=\frac{1}{2}$. The optimal value of $\lambda$ is determined
as the solution of $D'(\lambda)=0$, which yields $t$ as a solution od $f'(t)=\frac{f(b)-f(a)}{b-a}$. If $f$ does not have
a derivative everywhere in $(a,b)$, one can use methods relying on left and right derivatives. A discussion related to cases with non-uniform distribution
is out of scope of this paper.

}

\eitem

\qitem{$TH$ inequality with purely discrete measures}{Let $x_0<x_1<\ldots <x_n$, $n\geq 2$,  and let  $\mu(\{x_i\})=p_i$ where $p_i\in (0,1)$ and
$\sum p_i=1$, $p_i >0$.  The interval $[a,b]$ is here $[x_0,x_n]$. We can allow $t$ to be any point in the interval $(x_0,x_n)$; it can be one of points
$t_i$ with positive probability, or not.  Then (\ref{thh2}) with  a discrete measure $\mu$ reads:
%\lefteqn{

\beqn
\label{hadied}
\sum_{i=0}^{n} p_i f(x_i)\leq  & &  \frac{x_n-\sum_{i=0}^n p_ix_i }{x_n-x_0}f(x_0) + \frac{\sum_{i=0}^n p_ix_i-x_0}{x_n-x_0} f(x_n)  \\
                          & &  +\left(f(t)-\frac{x_n-t}{x_n-x_0}f(x_0) -\frac{t-x_0}{x_n-x_0}f(x_n)\right)\\
                          & & \times \left(\frac{1}{t-x_0}\sum_{x_i\leq t} p_i(x_i-x_0) + \frac{1}{x_n-t}\sum_{x_i >t} p_i (x_n-x_i)\right)
\eeqn

Although we will not discuss concrete examples, let us emphasize that all further results of this paper are also valid for discrete measures.

}

\eitem

\section{Quantifying  the  tightness via Karamata's theorem\label{kamata}}

As an introduction to the topic of this section, let us note  that all three inequalities that we considered so far are of the type

\be
\label{typin}
 \int f(x) \di G(x) \geq (\leq) \int f(x)\di  H(x)
 \ee
 where $G$ and $H$ are distribution functions of corresponding measures. Let $G$ be the distribution function that appears in
 the integral  $\int f(x) \di G(x)$ in inequalities $J$, $H$ and $TH$. The second measure is derived from $G$  as follows.

\smallskip

\qitem{Second measure in inequalities $J$, $H$ and $TH$}{\label{merin}
Let $c:=\int x\di G(x)$. The second measures are discrete  and derived  from $G$ as follows.

{\bf (J)}\  The second measure is the unit mass at $c$, with $H(x)=I_{[c,+\infty)}(x)$, and in these terms, the
Jensen's inequality can be written as $\int f(x) \di (G(x)-H(x))\geq 0$.

\smallskip

{\bf (H)}\  The second measure is concentrated at
points $a$ and $b$ with probabilities  $\frac{b-c}{b-a}$ and $\frac{c-a}{b-a}$ respectively, so
$H(x)=\frac{b-c}{b-a}I_{[a,b)}+I_{[b,+\infty)}$. This inequality is of the form $\int f(x) \di (G(x)-H(x))\leq 0$.

\smallskip

{\bf (TH)}  The second measure is  concentrated on the set $\{a,t,b\}$ with probabilities
$p_a$, $p_t$ and $p_b$ in  Corollary \ref{Colt}. The distribution function is
\[ H(x)= p_aI_{[a,+\infty)}(x) + p_tI_{[t,+\infty)}(x) +p_bI_{[b,+\infty)}(x), \]
and the inequality is of the form $\int f(x) \di (G(x)-H(x))\leq 0$.

}
\eitem

The next lemma gives some common properties of inequalities of type as in (\ref{typin}).

 \medskip

\qitem{Lemma}{\label{compp} Suppose that for measures $G$ and $H$  the inequality

\be
\label{typinp}
\int f(x) \di G(x) \geq \int f(x)\di  H(x)
\ee
holds with any convex function $f$ on $[a,b]$. Then

 \be
 \label{2cond}
 \int_{[a,b]} \di G(x)=\int_{[a,b]} \di H(x)\quad \mbox{and}\quad \int_{[a,b]} x\di G(x) = \int_{[a,b]} x\di H(x).
 \ee
Further, if $f$ is an affine function, the inequality (\ref{typin}) turns to equality.

\begin{proof}  The first equality follows upon plugging $f=1$ and $f=-1$ in (\ref{typinp}). For the second equality take $f(x)=x$ and $f(x)=-x$.
If  $f=\alpha x +\beta$, the statement above follows from (\ref{2cond}) using the linearity of integral.
\end{proof}

}

\eitem

  In the paper  \cite{karam32}, Jovan Karamata in the year 1932 presented  conditions for  two given measures so that
 the inequality (\ref{typinp}) holds with all  convex functions. This result is often wrongly attributed to Levin and
 Ste\v ckin  \cite{stelev1960}. In fact, \cite{stelev1960} was originally written by Ste\v ckin sixteen years after Karamata's paper,
  as  Supplement I in \cite{stec48},  with Karamata's paper \cite{karam32} in the list of references of \cite{stec48}. In several recently published papers,
  (for example \cite{rajb2014}), a related  result  is again rediscovered with  the name Ohlin's lemma,  after the paper \cite{ohlin1969} of the year 1969 in the context of
  application in actuarial area.

 \smallskip

\qitem{Theorem (Karamata \cite{karam32})}{\label{teka} Given two measures with distribution functions $G$ and $H$ and assuming conditions (\ref{2cond}),
the inequality (\ref{typinp})  holds for every  convex function $f$ if and only if for all  $u\in [a,b]$
\be
\label{iffa}
\fai (u):=\int_{[a,u]}(G(x)-H(x))\di x \geq 0\qquad \mbox{for all $u\in [a,b]$}
\ee
\hfill$\square$

}
\eitem
In the sequel we will refer to the function (\ref{iffa}) as   {\em Karamata's function}.
It is well known  (since as early as \cite{halipo29} and \cite{karam32}) that a function $f$ which is convex on $[a,b]$ can be uniformly approximated
by functions of the form
\be
\label{appconv}
x\mapsto \alpha x +\beta + \sum_{i=1}^n c_i |x-u_i|,\quad c_i>0,\ u_i\in (a,b)
\ee
 This observation together with the next lemma, completes the proof of Theorem \ref{teka}.

\qitem{Lemma}{\label{tekal} Under conditions (\ref{2cond}), the Karamata's function $u\mapsto \fai(u)$
can be represented as the residual in (\ref{typinp}) with the function $x\mapsto |x-u|$:
\be
\label{tekali}
\fai (u)=\int_{[a,b]} |x-u| \di G(x) -\int_{[a,b]} |x-u| \di H(x)
\ee

\begin{proof} Let $F(x):= G(x)-H(x)$. Then by conditions (\ref{2cond}) we have that $\int_{[a,b]}\di F(x)=\int_{[a,b]}x\di F(x)=0$, and

\beq
 \label{equiva}
  \int_{[a,b]} |x-u|\di F(x) &=& \int_{[a,u]} (u-x)\di F(x)+\int_{(u,b]} (x-u)\di F(x)\nonumber\\
                             & =& 2uF(u)  -2\int_{[a,u]}x\di F(x).
  \eeq

 Further, an integration by parts yields
 \[\int_{[a,u]}x\di F(x)=xF(x)\limitsi{a_{-}}{u} - \int_{[a,u]}F(x)\di x = uF(u)-\int_{[a,u]} F(x) \di x, \]
 so, from  (\ref{equiva}) it follows

 \[ \int_{[a,u]} F(x)\di x = \int_{[a, b]}|x-u| \di F(x),\]
which ends the proof.

\end{proof}

}
\eitem

\qitem{Average residuals}{In order to compare sharpness and tightness of two inequalities, we need to have a representative
measure for  the  size of residuals of an inequality itself, with no particular function attached.  In view of Lemma \ref{tekal},
a natural choice is the mean value of the Karamata's function $\fai$.  For  inequality $\cI$  which  satisfies  conditions of Karamata's
theorem,  we define the average residual as

\be
\label{mrs}
 AR (\cI) = \frac{1}{b-a}\int_a^b \fai(u)\di u
  \ee

 For comparing errors in two inequalities
 $\cI$ and $\cI_0$ on the same interval,  we  define relative average residual of $\cI$ with respect to $\cI_0$ as

 \be
 \label{rmrs}
 RAR (\cI, \cI_0)= \frac{AR (\cI)}{AR (\cI_0)}= \frac{\int_a^b \fai (u)\di u}{\int_a^b \fai_{0} (u)\di u}.
 \ee

For a concrete convex function, the size of residual  depends on the second derivative
(see explicit dependence formulae in \cite{repje} for some particular cases)
or some other measures of convexity. Although the residual and relative residual here can be calculated directly, a representation
of residuals in terms of  Karamata's function is meaningful to reveal to which extent the average residuals reflect particular ones.
The next theorem gives the relationship between residual (with given function) and  Karamata's function.

}

\eitem

\smallskip

\qitem{Theorem}{\label{rescf} Let $R(f,\cI)$ be the residual in inequality $\cI$, with given measures $G$ and $H$  and
with a twice differentiable  convex function $f$ on the interval $[a,b]$.  Then,

\be
\label{rfi}
R(f,\cI) = \frac{1}{2} \int_a^b f''(u) \fai(u) \di u= \frac{1}{2}\fai(\theta) (f'(b)-f'(a)),
\ee
for some $\theta\in (a,b)$.

\begin{proof} Let $h(x) = \int_a^b f''(u) |x-u| \di u$.
Performing the integration by parts on intervals $[a,u]$ and $[u,b]$ separately and adding, we find that $ h(x)= 2f(x) +g(x)$,
where $g(x)$ is affine. Therefore,

\be
\label{rfi1}
 f(x)=\frac{1}{2}h(x) - \frac{1}{2}g(x).
 \ee
Applying the inequality $\cI$ on both sides in (\ref{rfi1}), and using lemma \ref{tekal} and second statement in lemma \ref{compp}, we get

\be
\label{rexp}
R(f,\cI) =  \int_{[a,b]} f(x) \di (G(x)-H(x))=\frac{1}{2} \int_a^b f''(u) \fai(u) \di u.
\ee

Since $\fai$ is continuous and $f''\geq 0$, the second equality in (\ref{rfi}) follows from  an integral mean value theorem.
\end{proof}

For a given convex function $f$ in inequality $\cI$, we    define a relative residual with respect to $\cI_0$  as
\be
\label{rrf}
RR(f,\cI,\cI_0)= \frac{R(f,\cI)}{R(f,\cI_0)}=
\frac{\int_a^b f''(u) \fai(u) \di u}{\int_a^b f''(u) \fai_{0}(u) \di u}
\ee

}

\eitem

If the function $f$ is not twice differentiable, the following theorem gives a possibility of approximate residuals in the form as above.

\smallskip

\qitem{Theorem}{\cite[Theorem 2]{koliha2003}}{\label{koliht}   If $f$ is convex on $[a,b]$, then for any $\eps>0$ there
exists a convex $C^{\infty}$-function $\hat{f}$ such that $|f(x)-\hat{f}(x)|\leq \eps$ for all $x\in [a,b]$.\hfill$\square$

}

 \eitem

 Let  $\hat{f}=\hat{f}_{\eps}$ be an approximation for $f$ as in the theorem above, with some $\eps>0$.  It is not difficult to show that

\[ |R(f,\cI)-R(\hat{f}_{\eps},\cI)| \leq 2\eps\]
and
\[|RR (f,\cI,\cI_0) - RR (\hat{f}_{\eps},\cI,\cI_0)| \leq 2 \frac{R(\hat{f}_{\eps},\cI) + R(\hat{f}_{\eps},\cI_0)}{R(f,\cI_0)R(\hat{f}_{\eps},\cI_0)}\eps. \]

Therefore, formulae (\ref{rfi}) and (\ref{rrf})  can be used with $\hat{f}$ in place of $f$, with  small enough $\eps$ to achieve an arbitrary
small error of approximation.

\section{Numerical evidence: Graphic contents and tables}

In this section we compare the residuals of inequalities $J$, $H$ and $TH$, using the methodology  presented in the section \ref{kamata}.
The figures 1-3 are obtained by Maple calculation of Karamata's function in an equivalent form adopted for measures with densities:

\begin{itemize}
\item[(J)]\ $\dsp \fai(u)=\int_{[a,u]} (u-x)\di G(x) - (u-c) I_{[c,b]}(u)   \quad (c=\int x\di G(x) ),$
\item[(H)]\ $\dsp \fai(u)= \frac{u-a}{b-a} \int_{[a,b]} (b-x)\di G(x) - \int_{[a,u]} (u-x) \di G(x),$
\item[(TH)]\ $\dsp \fai(u) = \left( \frac{u-a}{t-a}\int_{[a,t]}(t-x)\di G(x)-\int_{[a,u]}(u-x)\di G(x)\right)I_{[a,t]}(u)$\\
          \rule{2.8em}{0em}     $\dsp+\left(\frac{u-t}{b-t}\int_{(t,b]}G(x)\di x -\int_{(t,u]}G(x)\di x \right)I_{(t,b]}\quad  (t=\frac{1}{2}).$
  \end{itemize}
Here $G$ is a main measure and $H$ is given in explicit form in terms of $G$, according to formulae in \ref {merin}.

In all cases we set $a=0$, $b=1$, and we consider three distributions:

\begin{itemize}
\item Uniform distribution on $[0,1]$, $G(x)=x$, $x\in [0,1]$.
\item Beta $(2,2)$ distribution, $G(x)= x^2(3-2x), x\in [0,1]$.
\item Exponential reduced to $[0,1]$ $G(x)=\frac{1-e^{-\lambda x}}{1-e^{-\lambda}}$, $\lambda=1$, $x\in [0,1]$.
\end{itemize}

Figures 1-3 show graphs of Karamata's functions $\fai$ for $J$-inequality (with spike), H-inequality (the largest)
and TH-inequality (lowest).

\begin{figure}[!htb]
   \begin{minipage}{0.33\textwidth}
     \centering
     \includegraphics[width=.8\linewidth]{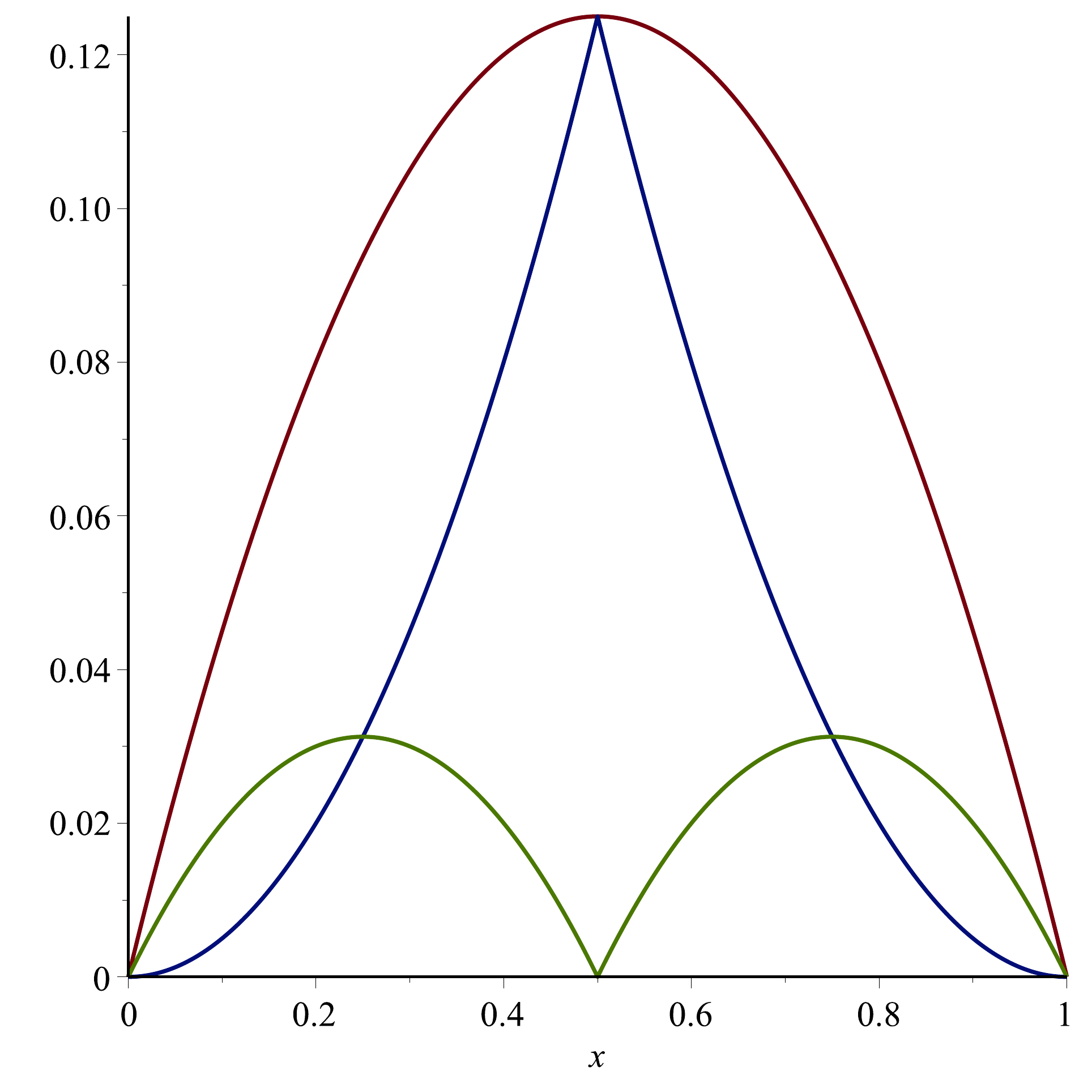}
     \caption{Uniform}\label{unifK}
   \end{minipage}\hfill
   \begin{minipage}{0.33\textwidth}
     \centering
     \includegraphics[width=.8\linewidth]{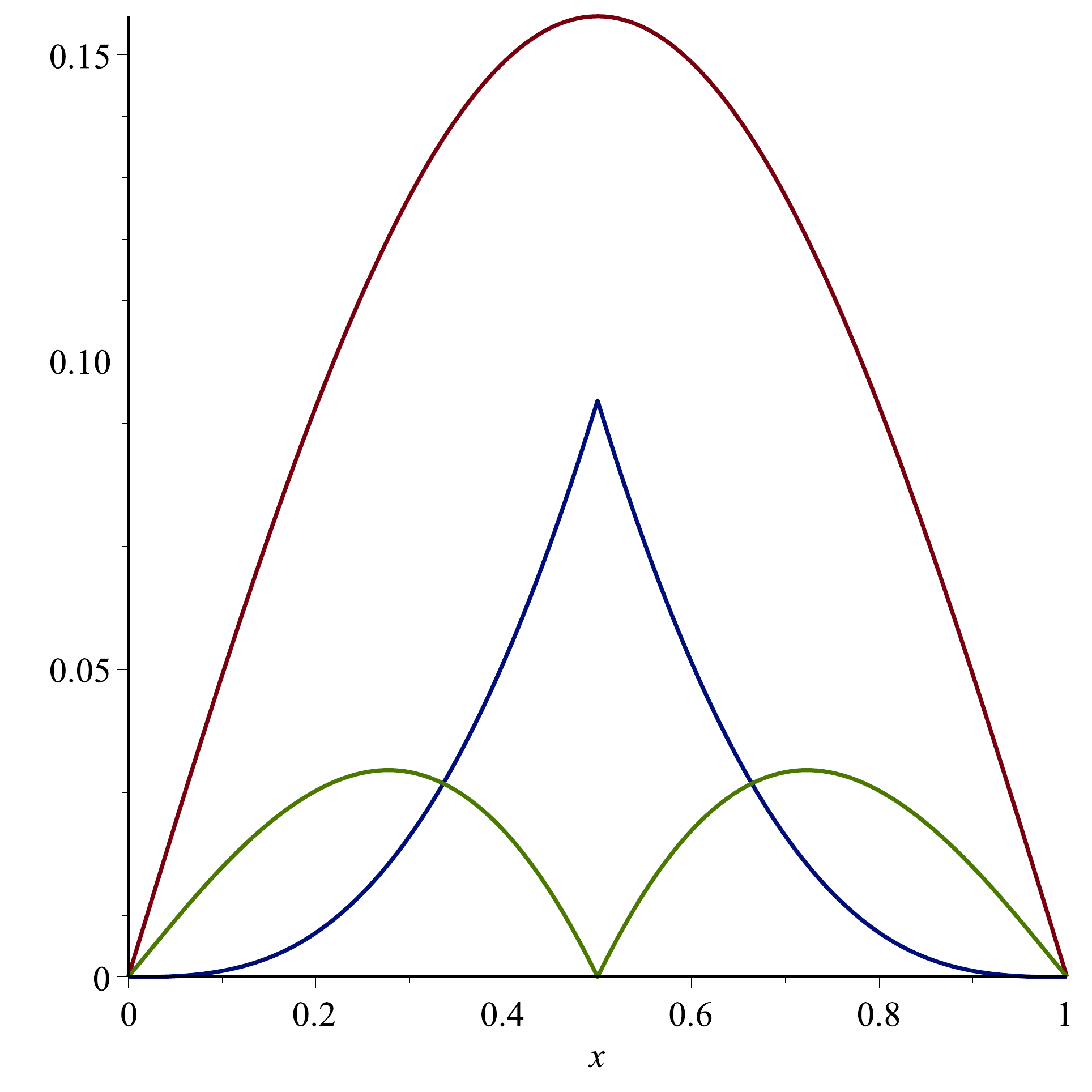}
     \caption{Beta}\label{betaK}
   \end{minipage}\hfill
   \begin{minipage}{0.33\textwidth}
     \centering
     \includegraphics[width=.8\linewidth]{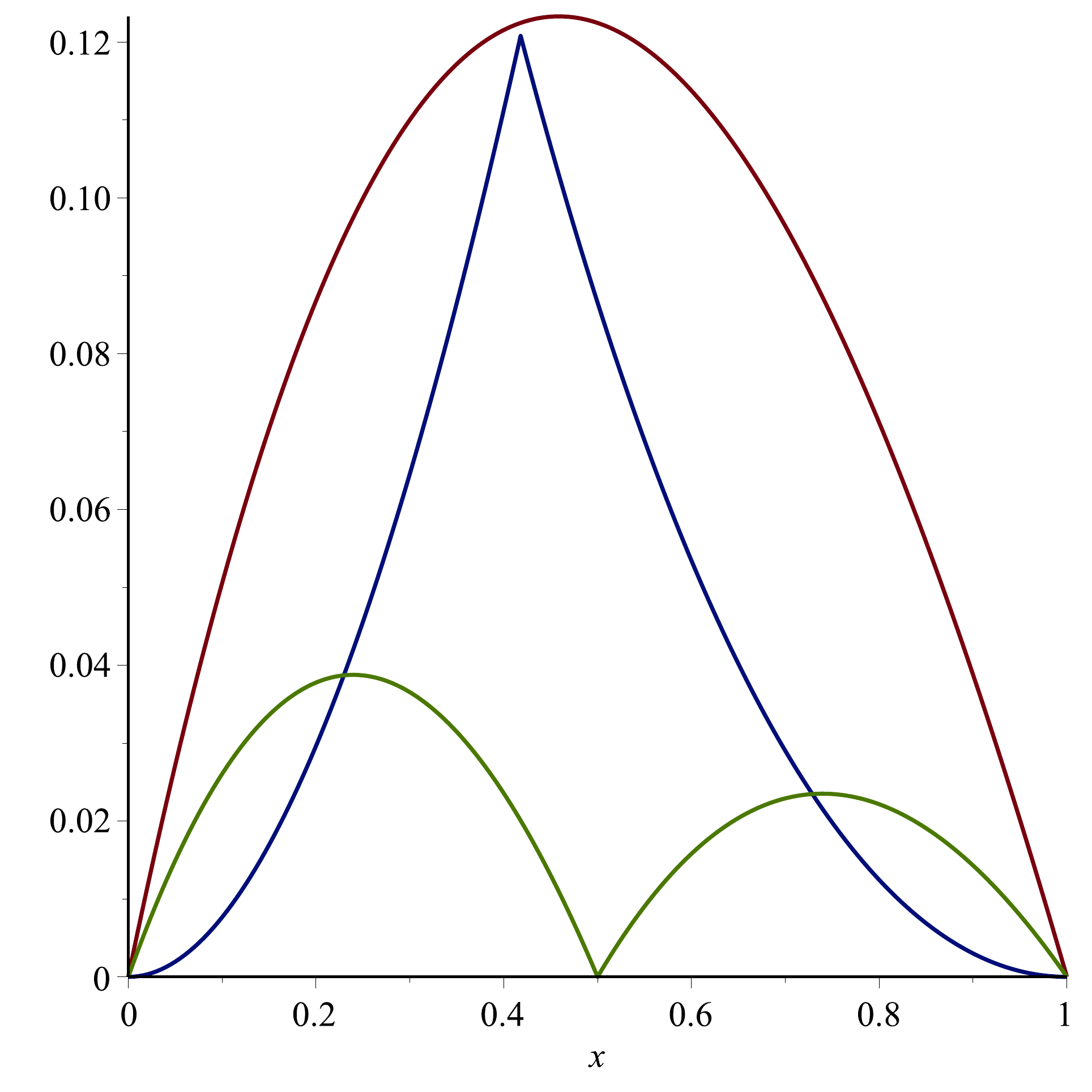}
     \caption{Exponential}\label{exprK}
   \end{minipage}
   \end{figure}

Since the domain  is the interval $[0,1]$, the average residual size is numerically
equal to the area between the $x$-axis and the graph. It is obvious that the area under TH curve is the smallest, in each od three examples with
different measures. This is confirmed in in the next table, where  we  present average residuals as in (\ref{mrs}).

\smallskip

\begin{center}
\begin{tabular}{l|ccc}\refstepcounter{table}
  &\multicolumn{3}{c}{Distribution (on $[0,1]$)}  \\
Inequality\rule{.1em}{0em}
& Uniform  & Exp $(1)$ & Beta $(2,2)$   \\ \hline
 Jensen   & 42  & 25 & 40  \\
Classical H  & 83    & 100 & 82 \\
Tight H & 21     & 22 & 21  \\ %\hline
\end{tabular}
\end{center}

\smallskip

\centerline{Table \thetable : The values of $AR \times 10^{3}$}

\smallskip

Relative average residuals can be derived from Table 1. For example,
if $\cI$ is  $TH$ with uniform distribution and $\cI_0$ is $H$ with $\expd (1)$ reduced to $[0,1]$, then $RAR(\cI,\cI_0)=0.21$.

We conclude that the theory in Sections 2 and 3, together with examples in this section,
show an absolute superiority of the new tight approximation to $\int f(x)\di G(x)$, compared to  classical Hermite-Hadamard bonds,
with $f$ being convex. The numerical evidences presented  in  the table above indicates that also it might be the case in comparison to
Jensen's lower bonds, which can be a topic of another research.

\medskip

\bibliography{convexity-hh}

\begin{thebibliography}{10}

\bibitem{allas2015}
{\sc Allasia, G.}
\newblock Connections between {H}ermite-{H}adamard inequalities and numerical
  integration of convex functions i.
\newblock {\em Bull. Allahabad Math. Soc. 30\/} (2015), 211--237.

\bibitem{bull78}
{\sc Bullen, P.}
\newblock Error estimates for some elementary quadrature rules.
\newblock {\em Univ. Beograd. Publ. Elektrotehn. Fak. Ser. Mat.Fiz. No.
  602-633\/} (1978), 97--103.

\bibitem{durp}
{\sc Durrett, R.}
\newblock {\em Probability: {T}heory and examples}.
\newblock Cambridge University Press, 2010.

\bibitem{guess2002}
{\sc Guessab, A., and Schmeisser, G.}
\newblock Sharp integral inequalities of the {H}ermite-{H}adamard type.
\newblock {\em J. Approx. Theory 115\/} (2002), 260--288.

\bibitem{guess2018}
{\sc Guessab, A., and Semisalov, B.}
\newblock A multivariate version of {H}ammer's inequality and its consequences
  in numerical integration.
\newblock {\em Results Math. 73\/} (2018), Art. 33, 37 pp.

\bibitem{hamm58}
{\sc Hammer, P.~C.}
\newblock The midpoint method of numerical integration.
\newblock {\em Math.Mag 31\/} (1958), 97--103.

\bibitem{halipo29}
{\sc Hardy, G.~H., Littlewood, J.~E., and P\'olya, G.}
\newblock Some simple inequalities satisfied by convex function.
\newblock {\em Messenger Math. 58\/} (1929), 145--152.

\bibitem{jen06}
{\sc Jensen, J. L. W.~V.}
\newblock Sur les fonctions convexes et les in\'egalit\'es entre les valeurs
  moyennes.
\newblock {\em Acta Math. 30\/} (1906), 175--193.

\bibitem{karam32}
{\sc Karamata, J.}
\newblock Sur une in\'egalit\'e relative aux fonctions convexes.
\newblock {\em Publ. Math. Univ. Belgrade 1\/} (1932), 145--148.

\bibitem{koliha2003}
{\sc Koliha, J.~J.}
\newblock Approximation of convex functions.
\newblock {\em Real Anal. Exchange 29\/} (2003), 465--471.

\bibitem{stelev1960}
{\sc Levin, V.~I., and Ste\v{c}kin, S.~B.}
\newblock Inequalities.
\newblock {\em Amer. Math. Soc. Transl 14\/} (1960), 1--22.

\bibitem{luppef76}
{\sc Lupa\c{s}, A.}
\newblock A generalization of {H}adamard inequalities for convex functions.
\newblock {\em Univ. Beograd. Publ. Elektrotehn. Fak. Ser. Mat. Fiz. 544-576\/}
  (1976), 115--121.

\bibitem{marolk}
{\sc Marshall, A.~W., Olkin, I., and Arnold, B.~C.}
\newblock {\em Inequalities: theory of majorization and its applications}.
\newblock Springer, 2009.
\newblock Springer series in {S}tatistics, second edition.

\bibitem{McShane37}
{\sc McShane, E.~J.}
\newblock Jensen's inequality.
\newblock {\em Bull. Amer. Math. Soc. 8\/} (1937), 521--527.

\bibitem{merc2009}
{\sc Mercer, P.~R.}
\newblock Hadamard's inequality and trapezoid rules for the
  {R}iemann-{S}tieltjes integral.
\newblock {\em J. Math. Anal. Appl. 344\/} (2008), 921--926.

\bibitem{mm99}
{\sc Merkle, M.}
\newblock Remarks on {O}strovski's and {H}adamard's inequality.
\newblock {\em Univ. Beograd. Publ. Elektrotehn. Fak. Ser. Mat. 10\/} (1999),
  113--117.

\bibitem{repje}
{\sc Merkle, M.}
\newblock Representation of the error term in {J}ensen's and some related
  inequalities with applications.
\newblock {\em J. Math. Analysis Appl. 231\/} (1999), 76--90.

\bibitem{mitlac}
{\sc Mitrinovi\'c, D.~S., and Lackovi\'c, I.}
\newblock Hermite and convexity.
\newblock {\em Aequations Math. 28\/} (1985), 229--232.

\bibitem{ohlin1969}
{\sc Ohlin, J.}
\newblock On a class of measures of dispersion with application to optimal
  reinsurance.
\newblock {\em ASTIN Bulletin 5\/} (1969), 249--266.

\bibitem{Olbs2015}
{\sc Olbry\'s, A., and Szostok, T.}
\newblock Inequalities of the {H}ermite-{H}adamard type involving numerical
  differentiation formulas.
\newblock {\em Results Math. 67\/} (2015), 403--416.

\bibitem{rajb2014}
{\sc Rajba, T.}
\newblock On the {O}hlin lemma for {Hermite-Hadamard-Fejer} type inequalities.
\newblock {\em Math. Inequal. Appl 17\/} (2014), 557--571.

\bibitem{stec48}
{\sc Ste\v{c}kin, S.~B.}
\newblock {\em Supplement 1: Inequalities for convex functions. In: {G. H.
  H}ardy, {J. E. L}ittlewood, {G. P}\'olya: {I}nequalities (in {R}ussian),
  translated from original by {V. I. L}evin, with supplements by {V. I.} Levin
  and {S. B. S}te\v ckin}.
\newblock Gosudarstvenoe izdatelstvo inostrannoi literaturi, Moskva, 1948,
  pp.~361--367.

\end{thebibliography}
\bibliographystyle{acm}

\bigskip

\noindent  Milan Merkle\\
University of Belgrade, School of Electrical Engineering\\
Bulevar kralja Aleksandra 73, 11120 Beograd, Serbia\\
emerkle@etf.rs

\medskip

\noindent Zoran D. Mitrovi\'c\\
University of Banja Luka, Faculty of Electrical Engineering\\
Patre 5, 78000 Banja Luka, Bosnia and Herzegovina\\
zoran.mitrovic@etf.unibl.org

\end{document}